\documentclass[a4paper,10pt]{amsart}
\usepackage{amsfonts}
\usepackage{float}
\usepackage[all]{xy}

\theoremstyle{plain}
\newtheorem{theorem}{Theorem}[section]
\newtheorem{prop}[theorem]{Proposition}
\newtheorem{lem}[theorem]{Lemma}
\newtheorem{corol}[theorem]{Corollary}

\newtheorem*{corolpresentationA}{Corollary \ref{corol:presentationA}}
\newtheorem*{theorempolyreg}{Theorem \ref{theorem:polyreg}}
\newtheorem*{corolregularinAQ}{Corollary \ref{corol:regularinAQ}}
\newtheorem*{theoremmultintubes}{Theorem \ref{theorem:multintubes}}
\newtheorem*{theoremmultinwild}{Theorem \ref{theorem:multinwild}}

\theoremstyle{definition}
\newtheorem{defi}[theorem]{Definition}

\newtheorem{rmq}[theorem]{Remark}
\newtheorem{exmp}[theorem]{Example}

\def\cone{{\rm{cone}\,}}
\def\coker{{\rm{coker}\,}}
\def\Hom{{\rm{Hom}}}
\def\Ext{{\rm{Ext}}}
\def\End{{\rm{End}}}
\def\Gr{{\rm{Gr}}}

\def\add{{\rm{add}}}

\def\dim{{\rm{dim}\,}}
\def\ddim{{\textbf{dim}\,}}
\def\rep{{\rm{rep}}}
\def\ind{{\rm{ind}}}

\def\Ob{{\rm{Ob}\,}}

\def\kQmod{{\emph{kQ}\textrm{-mod}}}

\def\Ai{{\mathbb{A}_{\infty}}}
\def\Aii{{\mathbb{A}_{\infty}^{\infty}}}

\def\A{{\mathbb{A}}}
\def\P{{\mathbb{P}}}
\def\N{{\mathbb{N}}}
\def\Z{{\mathbb{Z}}}
\def\C{{\mathbb{C}}}
\def\Q{{\mathbb{Q}}}

\def\CC{{\mathcal{C}}}

\def\<{\left<}
\def\>{\right>}
\def\ens#1{\left\{ #1 \right\}}
\def\fl{{\longrightarrow}\,}

\begin{document}
\title[Generalized Chebyshev Polynomials]{Cluster multiplication in regular components via generalized Chebyshev polynomials}
\date{}
\author{\textsc{G. Dupont}}
\address{Universit\'e de Lyon~; \\
Universit\'e Lyon 1~; \\
Institut Camille Jordan CNRS UMR 5208~; \\
43, boulevard du 11 novembre 1918\\
F-69622 Villeurbanne Cedex.}
\email{dupont@math.univ-lyon1.fr}

\begin{abstract}
	We introduce a multivariate generalization of normalized Chebyshev polynomials of the second kind. We prove that these polynomials arise in the context of cluster characters associated to Dynkin quivers of type $\mathbb A$ and representation-infinite quivers. This allows to obtain a simple combinatorial description of cluster algebras of type $\mathbb A$. We also provide explicit multiplication formulas for cluster characters associated to regular modules over the path algebra of any representation-infinite quiver.
\end{abstract}

\maketitle

\setcounter{tocdepth}{1}
\tableofcontents

\begin{section}{Introduction}
	Cluster algebras were introduced by Fomin and Zelevinsky in \cite{cluster1} in order to give a combinatorial framework for studying positivity in algebraic groups and canonical bases in quantum groups. Since, they were subject to various developments in areas like combinatorics, Lie theory, Teichm\"uller theory and quiver representations.
	
	By definition, a (simply-laced, coefficient-free) cluster algebra is a commutative algebra generated by a set of variables, called \emph{cluster variables}, consisting of possibly overlapping sets, all of the same finite cardinality, called \emph{clusters}. The initial data for constructing a cluster algebra is a pair $(Q,\textbf u)$ where $Q=(Q_0,Q_1)$ is a quiver without loops and 2-cycles, where $Q_0$ denotes the (finite) set of vertices, $Q_1$ the (finite) set of arrow and $\textbf u=(u_i, i \in Q_0)$ is a $Q_0$-tuple of indeterminates over $\Q$. The corresponding cluster algebra is denoted by $\mathcal A(Q)$. It is a subring of the ring $\Z[\textbf u^{\pm 1}]$ of Laurent polynomials in the $u_i$'s. If $Q$ has no oriented cycles, it is called \emph{acyclic} and $\mathcal A(Q)$ is called an \emph{acyclic cluster algebra}. 
	
	Categorifications of cluster algebras via representation theory of quivers were initiated in \cite{MRZ}. Since, it appeared that triangulated 2-Calabi-Yau categories provide a fruitful way for categorifying a wide class of cluster algebras (see \cite{Keller:categorification} and references therein). Given such a categorification of a cluster algebra $\mathcal A$, one has a notion of \emph{cluster character} in the sense of Palu, which allows to realize explicitly cluster variables in $\mathcal A$ from objects in the corresponding category \cite{Palu,FK}.
		
	The relevant triangulated 2-Calabi-Yau category for studying an acyclic cluster algebra $\mathcal A(Q)$ is the \emph{cluster category} $\CC_Q$ introduced in \cite{BMRRT} (see also \cite{CCS1}). The corresponding cluster character, in this case, is the so-called \emph{Caldero-Chapoton map} first introduced in \cite{CC} for Dynkin quivers and later generalized to the acyclic case in \cite{CK2}. The Caldero-Chapoton map is an explicit map $X_?:\Ob(\CC_Q) \fl \Z[\textbf u^{\pm 1}]$ from the set of objects in $\CC_Q$ to the ring of Laurent polynomials in $\textbf u$. Caldero and Keller proved that if $Q$ is any acyclic quiver, then cluster variables in $\mathcal A(Q)$ are precisely the characters $X_M$ where $M$ is an indecomposable rigid (that is, such that $\Ext^1_{\CC_Q}(M,M)=0$) object in $\CC_Q$ \cite{CK2}.

	Besides this property, the Caldero-Chapoton map also endows the acyclic cluster algebra $\mathcal A(Q)$ with a structure of Hall algebra on the cluster category $\CC_Q$. Namely, if $M$ and $N$ are two objects in the cluster category $\CC_Q$ such that $\Ext^1_{\CC_Q}(M,N) \neq 0$, the product $X_MX_N$ can be expressed as a $\Q$-linear combination of $X_Y$'s where $Y$ runs over isoclasses of middle terms of triangles involving $M$ and $N$. This property was first proved for Dynkin quivers in \cite{CK1}. It was generalized to arbitrary acyclic quivers when $\Ext^1_{\CC_Q}(M,N)$ is one-dimensional in \cite{CK2} and finally proved in general in \cite{XX}. Recently, Palu has obtained a similar multiplication formula for arbitrary cluster characters on Hom-finite 2-Calabi-Yau triangulated categories \cite{Palu:multiplication}.

	The very first illustration of the structure of Hall algebra of the cluster algebra was given by Caldero and Chapoton in terms of (finitely generated) modules over the path algebra $kQ$ of $Q$ where $Q$ is an acyclic quiver and $k$ is an algebraically closed field. The authors observed that if $M$ is a non-projective indecomposable $kQ$-module, then $X_MX_{\tau M}=X_B+1$ where $B$ is the middle term of the almost split sequence ending at $M$ \cite[Proposition 3.10]{CC}. When the considered module $M$ belongs to an homogeneous tube in the Auslander-Reiten quiver $\Gamma(kQ)$ of $kQ$-mod, it was observed by Caldero and Zelevinsky that the above multiplication property for almost split sequences  gives rise to normalized Chebyshev polynomials of the second kind \cite{CZ} (see Section \ref{section:poly} for definitions).
		
	In this paper, we prove that if $Q$ is any quiver of infinite representation type or if $Q$ is of Dynkin type $\A$ with a linear orientation, the multiplication property for almost split sequences gives rise to so-called \emph{cluster-mesh relations} and to a family of polynomials called \emph{generalized Chebyshev polynomials} which are multivariate generalizations of normalized Chebyshev polynomials of the second kind. These polynomials are also likely to have an interplay with the theory of orthogonal polynomials in several variables.

	The initial aim of this paper was to give a simplified version of Xiao-Xu's multiplication theorem for two characters $X_M$ and $X_N$ in the case where $M$ and $N$ are indecomposable modules in a same regular component of the Auslander-Reiten quiver $\Gamma(kQ)$ where $Q$ is any acyclic quiver. This is achieved using generalized Chebyshev polynomials. The statement we provide is less general but slightly more precise that the one given by Xiao and Xu in \cite{XX}.
		
	Studying these polynomials, it turned out that Generalized Chebyshev polynomials can also be used in order to give a simple presentation of a cluster algebra of Dynkin type $\A$. Also, using cluster-mesh relations, we introduce the notion of \emph{cluster-mesh algebra} of a stable translation quiver. These cluster-mesh algebras are of particular interest for studying certain subalgebras of an acyclic cluster algebra.
	
	Note that a deformed version of these generalized Chebyshev polynomials fitting to the context of acyclic cluster algebras with coefficients was recently introduced bu the author \cite{Dupont:qChebyshev}. 

	The paper is organized as follows. In Section \ref{section:background}, we give necessary background and state our main results. In Section \ref{section:poly}, we introduce cluster-mesh relations and generalized Chebyshev polynomials. In Section \ref{section:typeA}, we study cluster algebras associated to Dynkin quivers of type $\mathbb A$ through the lens of generalized Chebyshev polynomials. In Section \ref{section:regular}, we prove that generalized Chebyshev polynomials arise in the context of cluster characters associated to regular $kQ$-modules when $Q$ is of infinite representation type. In Section \ref{section:clustermeshalgebra}, we introduce the notion of cluster-mesh algebra of a translation quiver and investigate some particular examples. Finally, Section \ref{section:multiplication} is devoted to the proof of a multiplication formula for indecomposable modules in a regular component of the Auslander-Reiten quiver $\Gamma(kQ)$ of $kQ$-mod when $Q$ is any representation-infinite quiver. 
\end{section}

\begin{section}{Background and main results}\label{section:background}
	In this section, we briefly recall the necessary background concerning cluster categories and cluster characters. For a general overview of these notions, we refer to the survey \cite{Keller:categorification} and references therein.
	
	\begin{subsection}{Cluster categories and the Caldero-Chapoton map}
		Throughout the paper, $k$ denotes an algebraically closed field. Let $Q=(Q_0,Q_1)$ be an acyclic quiver. We denote by $kQ$-mod the category of finitely generated left $kQ$-modules over the path algebra $kQ$ of $Q$ which we will usually identify with the category $\rep(Q)$ of finite dimensional representations of $Q$ over $k$. For any vertex $i \in Q_0$, we denote by $S_i$ the simple module associated to vertex $i$, by $P_i$ its projective cover and by $I_i$ its injective hull. We denote by $\tau$ the Auslander-Reiten translation on $kQ$-mod and by $\<-,-\>$ the Euler form on $kQ$-mod given by
		$$\<M,N\>=\dim \Hom_{kQ}(M,N)-\dim \Ext^1_{kQ}(M,N).$$
		
		Since $kQ$-mod is hereditary, $\<-,-\>$ is well defined on the Grothendieck group $K_0(kQ)$ of $kQ$-mod. We denote by $\ddim:kQ\textrm{-mod} \fl \Z^{Q_0}$ the \emph{dimension vector} given by $\ddim M=(\dim \Hom_{kQ}(P_i,M))_{i \in Q_0}$ for any $kQ$-module $M$. Note that $\ddim$ induces an isomorphism of abelian groups $K_0(kQ) \simeq \Z^{Q_0}$ identifying the class of the simple module $S_i$ with the $i$-th vector $\alpha_i$ of the canonical basis of $\Z^{Q_0}$ for any $i \in Q_0$.
		
		Let $D^b(kQ)$ the bounded derived category of $Q$ equipped with shift functor $[1]$ and Auslander-Reiten translation $\tau$. The cluster category, as introduced in \cite{BMRRT}, is the orbit category $\CC_Q = D^b(kQ)/F$ where $F=\tau^{-1}[1]$ is an auto-equivalence of $D^b(kQ)$. Keller proved that the cluster category is a triangulated category and that the canonical functor $D^b(kQ) \fl \CC_Q$ is triangulated \cite{K}. The set $\ind(\CC_Q)$ of indecomposable objects in $\CC_Q$ can be described as follows~:
		$$\ind(\CC_Q)=\ind(\kQmod) \sqcup \ens{P_i[1] | i \in Q_0}.$$
		It is proved in \cite{BMRRT} that the cluster category is 2-Calabi-Yau, meaning that there is a bifunctorial duality 
		$$\Ext^1_{\CC_Q}(X,Y) \simeq D\Ext^1_{kQ}(Y,X)$$
		for any two objects $X,Y$ in $\CC_Q$ where $D=\Hom_k(-,k)$. Moreover, if $X$ and $Y$ are $kQ$-modules, then
		$$\Ext^1_{\CC_Q}(X,Y) \simeq \Ext^1_{kQ}(X,Y) \oplus D\Ext^1_{kQ}(Y,X).$$
		
		It follows from the above equality that a $kQ$-module is rigid in $kQ$-mod if and only if it is rigid in $\CC_Q$. An object $M$ is called \emph{basic} if any two distinct direct summands of $M$ are non-isomorphic. An object $T$ in $\CC_Q$ is called \emph{cluster-tilting} if $T$ is rigid, basic and if $\Ext^1_{\CC_Q}(T,X)=0$ implies $X \in \add(T)$.
		
		The cluster category $\CC_Q$ is known to provide a fruitful categorification of the acyclic cluster algebra $\mathcal A(Q)$ with initial seed $(Q,\textbf u)$ (see e.g. \cite{BMRRT,BMRT,CC,CK2}). In \cite{CC,CK2}, the authors introduced a map 
		$$X_?:\Ob(\CC_Q) \fl \Z[\textbf u^{\pm 1}]$$
		from the set of objects in the cluster category to the ring of Laurent polynomials $\Z[\textbf u^{\pm 1}]$ containing the cluster algebra $\mathcal A(Q)$. This map is a normalized generating series for Euler characteristics of certain algebraic varieties, called grassmannians of submodules. Given a $kQ$-module $M$ and a dimension vector $\textbf e \in \N^{Q_0}$, the \emph{grassmannian of submodules of $M$ of dimension $\textbf e$} (also called \emph{quiver grassmannian} in the literature) is the set
		$$\Gr_{\textbf e}(M)=\ens{N \subset M \textrm{ submodule such that }\ddim N=\textbf e}.$$
		It is a projective variety as it is a closed subset of the grassmannian of $k$-vector spaces. It is thus possible to consider its Euler characteristics (of the underlying topological space if $k=\C$ is the field of complex numbers or of $l$-adic cohomology if $k$ is arbitrary). We denote it by $\chi(\Gr_{\textbf e}(M))$.
		\begin{defi}
			The \emph{Caldero-Chapoton map} on $\CC_Q$ is the map 
			$$X_?:\Ob(\CC_Q)\fl \Z[\textbf u^{\pm 1}]$$
			defined as follows~:
			\begin{enumerate}
				\item If $M$ is an indecomposable $kQ$-module, then
					\begin{equation}\label{eq:CCmap}
						X_M = \sum_{\textbf e} \chi(\Gr_{\textbf e}(M)) \prod_{i \in Q_0} u_i^{-\<\textbf e, \alpha_i\>-\<\alpha_i, \ddim M - \textbf e\>}~;
					\end{equation}
				\item If $M=P_i[1]$ for some $i \in Q_0$, then $$X_M=u_i~;$$
				\item For any two objects $M,N$ of $\CC_Q$, 
					$$X_{M \oplus N}=X_M X_N.$$
			\end{enumerate}
		\end{defi}
		 Note that $X_?$ is constant on isoclasses and equality (\ref{eq:CCmap}) also holds for decomposable $kQ$-modules. This map is now known to be a particular case of cluster character in the sense of Palu \cite{Palu}. For this reason, for any object $M$ in the cluster category, $X_M$ will sometimes be called the \emph{cluster character} or simply \emph{character} associated to $M$.

		The main motivation for studying the Caldero-Chapoton map is the following theorem~:
		\begin{theorem}[\cite{CK2}]\label{theorem:correspondenceCK2}
			Let $Q$ be an acyclic quiver. Then $X_?$ induces a 1-1 correspondence from the set of indecomposable rigid objects in $\CC_Q$ to the set of cluster variables in $\mathcal A(Q)$. Moreover, $X_?$ induces a 1-1 correspondence from the set of cluster-tilting objects in $\CC_Q$ to the set of clusters in $\mathcal A(Q)$.
		\end{theorem}
		
		In \cite{BMRT}, the authors provided a surjective map $\alpha$ from the set of cluster variables in $\mathcal A(Q)$ to the set of indecomposable rigid objects in $\CC_Q$. It is known that $X_?$ is a left inverse to the map $\alpha$. Thus, when one only consider rigid objects, both approaches developed respectively in \cite{CC,CK1,CK2} and in \cite{BMRT} are equivalent. Nevertheless, in our context, we will mainly deal with cluster characters associated to indecomposable regular modules over the path algebra of a representation-infinite quiver and it is known that all but finitely many such modules are non-rigid. From this point of view, the fact that the Caldero-Chapoton map allows to consider non-rigid objects in $\CC_Q$ is essential. Note that it will appear in Corollary \ref{corol:regularinAQ} that these characters associated to non-rigid objects may also belong to the cluster algebra.
		
		As we already mentioned, the Caldero-Chapoton map also endows $\mathcal A(Q)$ with a structure of Hall algebra on the cluster category. The following theorem, illustrating this fact, will be of particular interest throughout the paper~:
		\begin{theorem}[\cite{CK2}]\label{theorem:multiplicationCK2}
			Let $Q$ be an acyclic quiver and $M,N$ be two indecomposable objects in $\CC_Q$ such that $\Ext^1_{\CC_Q}(M,N) \simeq k$, then $X_MX_N=X_B+X_{B'}$ where $B$ and $B'$ are the unique (up to isomorphism) objects in $\CC_Q$ such that there exists non-split triangles 
			$$M \fl B \fl N \fl M[1] \textrm{ and } N \fl B'\fl M \fl N[1].$$
		\end{theorem}

		Actually, in this paper, we will mainly use the following classical module-theoretic interpretation of the above theorem which first appeared in \cite{Hubery:cluster}. We give a proof for completeness.
		\begin{corol}\label{corol:multiplicationCK2module}
			Let $Q$ be an acyclic quiver and $M,N$ be two indecomposable $kQ$-modules such that $\Ext^1_{kQ}(M,N) \simeq k$ and $\Ext^1_{kQ}(N,M)=0$. Then, $X_MX_N=X_E+X_B$ where $E$ is the unique (up to isomorphism) $kQ$-module such that there exists a non-split exact sequence  $$0 \fl N \fl E \fl M \fl 0$$ in $kQ$-mod and $B = \ker f \oplus \tau^{-1} \coker f$ where $f$ is any non-zero morphism in $\Hom_{kQ}(N,\tau M)$.
		\end{corol}
		\begin{proof}
			Let $M$ and $N$ be two indecomposable modules such that $\Ext^1_{kQ}(M,N) \simeq k$ and $\Ext^1_{kQ}(N,M)=0$. There exists a unique (up to isomorphism) $kQ$-module $E$ such that there is a non-split exact sequence of $kQ$-modules 
			$$0 \fl N \fl E \fl M \fl 0.$$
			This short exact sequence induces a non-split triangle
			$$N \fl E \fl M \fl N[1]$$
			in $D^b(kQ)$ which induces a non-split triangle 
			$$N \fl E \fl M \fl N[1]$$
			in $\CC_Q$.
			According to the Auslander-Reiten formula, there are isomorphisms of $k$-vector spaces $\Hom_{kQ}(N,\tau M)\simeq \Ext^1_{kQ}(M,N) \simeq k$. Fix a non-zero morphism $f \in \Hom_{kQ}(N,\tau M)$, then there is a non-split triangle 
			$$\tau M[-1] \fl N \fl \cone(f)[-1] \xrightarrow{f} \tau M$$
			in $D^b(kQ)$ where $\cone(f)=\ker f[1] \oplus \coker f$.
			Thus, there is a non-split triangle 
			$$M \fl N \fl \ker f \oplus \tau^{-1} \coker f \xrightarrow{f} M[1]$$
			in $\CC_Q$. Since, $\Ext^1_{\CC_Q}(M,N) \simeq \Ext^1_{kQ}(M,N) \oplus D\Ext^1_{kQ}(N,M) \simeq k$, Theorem \ref{theorem:multiplicationCK2} gives $X_MX_N=X_E+X_B$, which proves the corollary.
		\end{proof}
	\end{subsection}

	\begin{subsection}{Auslander-Reiten theory of representation-infinite quivers}
		We fix a representation-infinite acyclic quiver $Q$. We shall briefly recall classical results on the Auslander-Reiten theory of $kQ$. The Auslander-Reiten quiver of $kQ$-mod is denoted by $\Gamma(kQ)$. It contains infinitely many connected components. There is one component containing all the projective (resp. injective) $kQ$-modules which is called the \emph{preprojective} component (resp. \emph{preinjective} component). The components which are not preprojective are called \emph{regular} components. An indecomposable module is called \emph{regular} if it belongs to a regular component. A decomposable module is called \emph{regular} if all its indecomposable direct summands are regular. 
		
		The structure of regular components in $\Gamma(kQ)$ is well known. For a quiver $R=(R_0,R_1)$, we recall that the \emph{repetition quiver} $\Z R$ associated to $R$ is the quiver with vertex set $\Z \times R_0$ and whose arrows are given by $(n,i) \fl (n,j)$ and $(n,j) \fl (n+1,i)$ for any arrow $i \fl j$ in $R_1$. Let $\Ai$ be the quiver with vertex set $\N$ and with arrows $i \fl i+1$ for any $i \in \N$. Then $\Z\Ai$ is a stable translation quiver for the translation $\tau$ given by $\tau (i,n)=(i-1,n)$. Then, for every regular component $\mathcal R$ in $\Gamma(kQ)$ there exists an integer $p \geq 0$ such that $\mathcal R$ is of the form $\Z\Ai/(\tau^p)$ which is also denoted $\Z\Ai/(p)$ \cite[Section VIII.4]{ARS}. 
		
		If $p >0$, $\mathcal R$ is called a \emph{tube of rank $p$}. If $p=1$, $\mathcal R$ is called \emph{homogeneous} and if $p>1$, $\mathcal R$ is called \emph{exceptional}. If $Q$ is an affine quiver, then the regular components form a family of tubes parametrized by $\P^1(k)$ and at most three of these tubes are exceptional \cite{ringel:1099}. Moreover, it is known that if $\mathcal R$ is a tube in $kQ$-mod, then $\mathcal R$ is standard, that is, the full subcategory generated by indecomposable objects in $\mathcal R$ is equivalent to the mesh category of $\mathcal R$. If $Q$ is a wild quiver, then regular components in $\Gamma(kQ)$ are all of the form $\Z\Ai$ \cite{ringel:wild}.

		We now return to the case where $Q$ is any acyclic representation-infinite quiver and we fix a regular component $\mathcal R$ in $\Gamma(kQ)$ which of the form $\Z\Ai/(p)$ for some $p \geq 0$. An indecomposable $kQ$-module $X$ in $\mathcal R$ is called quasi-simple if it is not a direct summand of the middle term of an almost split sequence in $kQ$-mod. We denote by $R_i,i \in \Z$ the quasi-simple modules in $\mathcal R$ ordered such that $\tau R_i \simeq R_{i-1}$ and $R_{i +p } = R_i$ for any $i \in \Z$. For any indecomposable $kQ$-module $M$ in $\mathcal R$, there exists an unique sequence of monomorphisms $0 =R_i^{(0)} \fl R_i = R_i^{(1)} \fl R_i^{(2)} \fl \cdots \fl R_i^{(n)}=M$
		such that $R_i^{(k)}/R_i^{(k-1)}$ is quasi-simple in $\mathcal R$ for any $k=1, \ldots, n$. The quotients $R_i^{(k)}/R_i^{(k-1)}$ for $k=1, \ldots, n$ are called the \emph{quasi-composition factors} of $M$, $R_i$ is called the \emph{quasi-socle} of $M$ and $n$ is called the \emph{quasi-length} of $M$. Note that if we set $M_k=R_i^{(k)}/R_i^{(k-1)}$ for any $k=1, \ldots, n$, we have $\tau M_k=M_{k-1}$ for any $k=2, \ldots, n$.
	\end{subsection}
	
	\begin{subsection}{Main results}
		We will introduce the family of \emph{generalized Chebyshev polynomials} $\ens{P_n}_{n \geq 0}$ given by $P_0=1$, $P_1(x)=x$ and for any $n \geq 1$, $P_n$ is the polynomial in $n$ variables given by
		$$P_{n+1}(x_1, \ldots, x_{n+1})=x_{n+1}P_n(x_1, \ldots, x_n)-P_{n-1}(x_1, \ldots, x_{n-1}).$$
		We will prove in Lemma \ref{lem:3termes} that these polynomials are characterized by the \emph{cluster mesh relations}~:
		$$P_{n+1}(x_1,\ldots, x_{n+1})=\frac{P_n(x_1, \ldots, x_n)P_n(x_2,\ldots, x_{n+1})-1}{P_{n-1}(x_2, \ldots, x_n)}$$
		for any $n \geq 1$. 
		
		These generalized Chebyshev polynomials provide a simple presentation of a cluster algebra of Dynkin type $\A$ in the sense of \cite{cluster2}~:
		\begin{corolpresentationA}
			Let $r \geq 0$ and $\mathcal A$ be a cluster algebra of Dynkin type $\A_r$. Then there is an isomorphism of $\Z$-algebras 
			$$\mathcal A \simeq \Z[t_0, \ldots, t_r]/(P_{r+1}(t_0, \ldots, t_r)-1)$$
			where $t_0, \ldots, t_r$ are indeterminates over $\Z$.
		\end{corolpresentationA}
		
		We also prove that for any representation-infinite quiver $Q$ and any indecomposable regular $kQ$-module $M$, the character $X_M$ can be expressed as a generalized Chebyshev polynomials in the characters associated  to its quasi-composition factors. Namely, we prove~:
		\begin{theorempolyreg}
			Let $Q$ be a quiver of infinite representation type. Let $\mathcal R$ be a regular component of the form $\Z\Ai/(p)$ for some $p \geq 0$ in $\Gamma(kQ)$. Let $\ens{R_i|i \in \Z/p\Z}$ denote the set of quasi-simple modules in $\mathcal R$ ordered such that $\tau R_i \simeq R_{i-1}$ for any $i \in \Z/p\Z$.. Then, for any $n \geq 1$ and any $i \in \Z/p\Z$, we have
			$$X_{R_i^{(n)}}=P_n(X_{R_i}, \ldots, X_{R_{i+n-1}}).$$
		\end{theorempolyreg}
		
		As a corollary, we can prove that a certain family of cluster characters associated to non-rigid indecomposable objects does belong to the cluster algebra $\mathcal A(Q)$~:
		\begin{corolregularinAQ}
			Let $Q$ be an acyclic quiver and $M$ be an indecomposable regular $kQ$-module with rigid quasi-composition factors. Then, $X_M \in \mathcal A(Q)$.
		\end{corolregularinAQ}
		
		Combinatorics of generalized Chebyshev polynomials together with Theorem \ref{theorem:polyreg} allow to prove multiplication formulas for cluster characters associated to regular $kQ$-modules. Namely, we prove for affine quivers~:
		\begin{theoremmultintubes}
			Let $Q$ be an affine quiver and $\mathcal T$ be a tube of rank $p$ in $\Gamma(kQ)$. Let $R_i, i \in \Z$ denote the quasi-simple modules in $\mathcal T$ ordered such that $\tau R_i \simeq R_{i-1}$ and $R_{i+p} \simeq R_i$ for any $i \in \Z$. Let $m,n>0$ be integers and $j \in [0,p-1]$. Then, for every $k \in \Z$ such that $0 < j+kp  < n$ and $m > n-j-kp$, we have the following identity~:
			$$X_{R_{j}^{(m)}}X_{R_0^{(n)}}=X_{R_0^{(m+j+kp)}}X_{R_{j}^{(n-j-kp)}}+X_{R_0^{(j+kp-1)}}X_{R_{n+1}^{(m+j+kp-n-1)}}.$$
		\end{theoremmultintubes}
		and prove in Corollary \ref{corol:Hallproduct} that this can be interpreted as a Hall product.
		
		We also obtain a similar description for wild quivers~:
		\begin{theoremmultinwild}
			Let $Q$ be a wild quiver and $\mathcal R$ be a regular component in $\Gamma(kQ)$. Let $R_i, i \in \Z$ denote the quasi-simple modules in $\mathcal T$ ordered such that $\tau R_i \simeq R_{i-1}$ for any $i \in \Z$. Let $m,n>0$ be integers and $j \geq 0$ such that $0 < j < n$ and $m > n -j$. Then, we have the following identity~:
			$$X_{R_{j}^{(m)}}X_{R_0^{(n)}}=X_{R_0^{(m+j)}}X_{R_{j}^{(n-j)}}+X_{R_0^{(j-1)}}X_{R_{n+1}^{(m+j-n-1)}}.$$
		\end{theoremmultinwild}
		
	\end{subsection}
\end{section}

\begin{section}{Cluster-mesh relations and generalized Chebyshev polynomials}\label{section:poly}
	\begin{subsection}{Generalized Chebyshev polynomials of infinite rank}
		Throughout the paper, $\ens{t_i|i \in \Z}$ denotes a family of indeterminates over $\Z$. We also fix $\ens{x_{i}|i \in \Z}$ a family of indeterminates over $\Z$. We define a family $\ens{x_{i,n}|i \in \Z, n\geq 1}$ of elements in $\Q(x_{i}|i \in \Z)$ by setting, for any $n \geq 1$,
		\begin{equation}\label{eq:clustermeshrelation}
			x_{i,n}x_{i+1,n}=x_{i,n+1}x_{i+1,n-1}+1
		\end{equation}
		where $x_{i,0}=1$ and $x_{i,1}=x_i$ for any $i \in \Z$.
		
		We label vertices in $\Z\Ai$ with $\ens{x_{i,n}|i \in \Z, n\geq 1}$ by identifying $(i,n)$ with $x_{i,n}$ for any $i \in \Z$ and $n \geq 1$. Equality (\ref{eq:clustermeshrelation}) is called a \emph{cluster mesh relation}. The following figure enlightens the chosen terminology.
		
		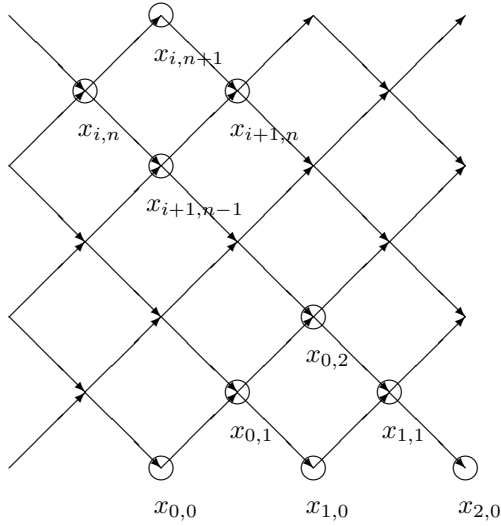
\begin{figure}[H]
			\setlength{\unitlength}{1mm}
			\begin{picture}(100,80)(-35,-10)
				\multiput(0,0)(20,0){3}{\multiput(0,0)(0,20){3}{\vector(1,1){10}}}
				\multiput(10,10)(20,0){3}{\multiput(0,0)(0,20){3}{\vector(1,-1){10}}}
				\multiput(0,0)(20,0){3}{\multiput(0,20)(0,20){3}{\vector(1,-1){10}}}
				\multiput(10,10)(20,0){3}{\multiput(0,0)(0,20){3}{\vector(1,1){10}}}
				
				\put(20,0){\circle{3}}
				\put(19,-6){$x_{0,0}$}
				
				\put(30,10){\circle{3}}
				\put(29,4){$x_{0,1}$}
				
				\put(40,20){\circle{3}}
				\put(39,14){$x_{0,2}$}
				
				\put(40,0){\circle{3}}
				\put(39,-6){$x_{1,0}$}
				
				\put(50,10){\circle{3}}
				\put(49,4){$x_{1,1}$}
				
				\put(60,0){\circle{3}}
				\put(59,-6){$x_{2,0}$}
				
				\put(20,60){\circle{3}}
				\put(19,54){$x_{i,n+1}$}
				
				\put(10,50){\circle{3}}
				\put(9,44){$x_{i,n}$}
				
				\put(30,50){\circle{3}}
				\put(29,44){$x_{i+1,n}$}
				
				\put(20,40){\circle{3}}
				\put(18,34){$x_{i+1,n-1}$}
			\end{picture}\label{figure:clustermeshrelation}
			\caption{Cluster-mesh relations}
		\end{figure}
		
		\begin{lem}\label{lem:rationalfunction}
			For any $n \geq 1$, there exists a rational function $P_n \in \Q(t_1, \ldots, t_n)$ such that
			$$x_{i,n}=P_n(x_i, \ldots, x_{i+n-1})$$
			for any $i \in \Z$.
		\end{lem}
		\begin{proof}
			We prove it by induction on $n$. By convention, we set $P_0=1$. If $n=1$, the result clearly holds and $P_1(t)=t$. Since for any $n \geq 1$ and for any $i \in \Z$, we have
			$$x_{i,n+1}=\frac{x_{i,n}x_{i+1,n}-1}{x_{i+1,n-1}}=\frac{P_n(x_i, \ldots, x_{i+n-1})P_n(x_{i+1}, \ldots, x_{i+n})-1}{P_{n-1}(x_{i+1}, \ldots, x_{i+n-1})}.$$
			Thus, $x_{i,n+1}=P_{n+1}(x_i, \ldots, x_{i+n})$ where 
			$$P_{n+1}(t_0, \ldots, t_n)=\frac{P_n(t_0, \ldots, t_{n-1})P_n(t_1, \ldots, t_n)-1}{P_{n-1}(t_1, \ldots, t_{n-1})}.$$
			This proves the lemma.
		\end{proof}

		We now prove that the rational functions $P_n$ satisfy a relation analogous to three terms recurrence relations in the theory of orthogonal polynomials in several variables (see e.g. \cite{Dunkl}).
		\begin{lem}\label{lem:3termes}
			For any $n \geq 2$, we have 
			$$P_n(t_0, \ldots, t_{n-1})=t_{n-1} P_{n-1}(t_0, \ldots, t_{n-2})-P_{n-2}(t_0, \ldots, t_{n-3})$$
		\end{lem}
		\begin{proof}
			We prove it by induction on $n$. For $n=2$, we compute directly
			$$P_2(t_0,t_1)=t_0t_1-1=t_1P_1(t_0)-P_0.$$
			According to the proof of lemma \ref{lem:rationalfunction}, we know that
			$$P_{n+1}(t_0,\ldots, t_n)= \frac{P_n(t_0, \ldots, t_{n-1})P_n(t_1, \ldots, t_n)-1}{P_{n-1}(t_1, \ldots, t_{n-1})}.$$
			By induction, $P_n(t_1, \ldots, t_n)=t_nP_{n-1}(t_1, \ldots, t_{n-1})-P_{n-2}(t_1, \ldots, t_{n-2})$ so that we get
			$$P_{n+1}(t_0,\ldots, t_n)=t_nP_n(t_0, \ldots, t_{n-1})-\frac{P_n(t_0, \ldots, t_{n-1})P_{n-2}(t_1, \ldots, t_{n-2})+1}{P_{n-1}(t_1, \ldots, t_{n-1})}.$$
			Since 
			$$P_n(t_0, \ldots, t_{n-1})P_{n-2}(t_1, \ldots, t_{n-2})+1=P_{n-1}(t_0, \ldots, t_{n-2})P_{n-1}(t_1, \ldots, t_{n-1}),$$
			we get
			$$P_{n+1}(t_0, \ldots, t_n)=t_n P_n(t_0, \ldots, t_{n-1}) - P_{n-1}(t_0, \ldots, t_{n-2})$$
			which proves the lemma.
		\end{proof}

		As a corollary, we obtain the polynomiality of the $P_n$ and we can moreover provide an explicit formula for these polynomials using determinental expressions.
		\begin{corol}
			For any $n \geq 1$, $P_n(t_0, \ldots, t_n)$ is the polynomial with integral coefficients given by 
			$$P_n(t_0, \ldots, t_n)=\det \left[\begin{array}{ccccccc}
			      t_n & 1 &&& (0)\\
			      1 & \ddots & \ddots \\
			      & \ddots & \ddots & \ddots \\
			      & & \ddots & \ddots & 1 \\
			      (0)& & & 1 & t_0
			\end{array}\right].$$
		\end{corol}
		\begin{proof}
			The identity is obtained by induction, noticing that the expansion of the determinent with respect to the first column is precisely the relation given in Lemma \ref{lem:3termes}. Polynomiality follows obviously from this identity. 
		\end{proof}
		
		\begin{defi}
			For any $n \geq 1$, $P_n$ is called the \emph{$n$-th generalized Chebyshev polynomial of infinite rank} of simply \emph{$n$-th generalized Chebyshev polynomial}.
		\end{defi}
		
		\begin{exmp}
			The first five generalized Chebyshev polynomials of infinite rank are
			$$\begin{array}{|c|c|}
				\hline 
				P_1(t_0) & t_0 \\
				\hline 
				P_2(t_0,t_1) & t_0t_1-1 \\
				\hline 
				P_3(t_0,t_1,t_2) & t_0 t_1 t_2 - t_0 - t_2 \\
				\hline 
				P_4(t_0,t_1,t_2,t_3) & t_0 t_1 t_2 t_3 - t_0 t_1 - t_2 t_3 - t_0 t_3 + 1 \\
				\hline 
				P_5(t_0,t_1,t_2,t_3,t_4) & t_0 t_1 t_2 t_3 t_4 - t_0 t_1 t_2 - t_0 t_1 t_4 - t_2 t_3 t_4 + t_2 - t_0 t_3 t_4 + t_0 + t_4 \\
				\hline 
			\end{array}$$
		\end{exmp}
	\end{subsection}

	\begin{subsection}{Generalized Chebyshev polynomials of finite ranks}
		For any $p \geq 1$, the abelian group $p\Z$ acts by homomorphisms of $\Z$-algebras on $\Z[x_{i}|i \in \Z]$ by 
		$$kp.x_{i}=x_{i+kp}$$
		for any $i,k \in \Z$. We denote by 
		$$\pi_p:\Z[x_{i}|i \in \Z] \fl \Z[x_{i}|i \in \Z]/p\Z$$
		the canonical map.
		
		Then, there exists a unique polynomial $P_{n,p} \in \Z(t_0, \ldots, t_{p-1})$ such that
		$$\pi_p(x_{i,n})=P_{n,p}(\pi_p(x_{i}), \ldots, \pi_p(x_{i+n-1}))$$
		for any $i \in \Z$, $n \geq 1$. 
		
		\begin{defi}
			For any $n \geq 1$ and $p \geq 1$, $P_{n,p}$ is called the \emph{$n$-th generalized Chebyshev polynomial of rank $p$}.
		\end{defi}
		
		We recall that for any $n \geq 1$, the classical \emph{$n$-th normalized Chebyshev polynomial of the second kind} is the polynomial $S_n(x) \in \Z[x]$ characterized by $S_0=1$, $S_1(x)=x$ and $S_{n+1}(x)=xS_n(x)-S_{n-1}(x)$ for any $n \geq 1$. The following straightforward lemma justifies the terminology of generalized Chebyshev polynomials~:
		\begin{lem}
			For any $n \geq 0$, the $n$-th generalized Chebyshev polynomial of rank 1 is the $n$-th normalized Chebyshev polynomial of the second kind.
		\end{lem}
		
		\begin{exmp}
			The first five generalized Chebyshev polynomials of rank one are 
			$$\begin{array}{|c|c|}
				\hline 
				P_{1, 1}(t_0) & t_0 \\
				\hline 
				P_{2, 1}(t_0) & t_0^2-1 \\
				\hline 
				P_{3, 1}(t_0) & t_0^3 - 2t_0 \\
				\hline 
				P_{4, 1}(t_0) & t_0^4- 3t_0^2+ 1 \\
				\hline 
				P_{5, 1}(t_0) & t_0^5- 4t_0^33t_0\\
				\hline 
			\end{array}$$
	
			The first five generalized Chebyshev polynomials of rank two are 
			$$\begin{array}{|c|c|}
				\hline 
				P_{1, 2}(t_0,t_1) & t_0 \\
				\hline 
				P_{2, 2}(t_0,t_1) & t_0 t_1 - 1 \\
				\hline 
				P_{3, 2}(t_0,t_1) & t_0^2 t_1 - 2t_0 \\
				\hline 
				P_{4, 2}(t_0,t_1) &    t_0^2  t_1^2  - 3 t_0 t_1 + 1 \\
				\hline 
				P_{5, 2}(t_0,t_1) & t_0^3  t_1^2  - 4 t_0^2 t_1 + 3t_0\\
				\hline 
			\end{array}$$
			
			The first five generalized Chebyshev polynomials of rank three are
			$$\begin{array}{|c|c|}
				\hline 
				P_{1, 3}(t_0,t_1,t_2) & t_0 \\
				\hline 
				P_{2, 3}(t_0,t_1,t_2) & t_0 t_1 - 1 \\
				\hline 
				P_{3, 3}(t_0,t_1,t_2) & t_0 t_1 t_2 - t_0 - t_2 \\
				\hline 
				P_{4, 3}(t_0,t_1,t_2) &      t_2 t_0^2  t_1 - t_0 t_1 - t_2 t_0 - t_0^2  + 1 \\
				\hline 
				P_{5, 3}(t_0,t_1,t_2) & t_0^2  t_1^2  t_2 - t_0^2  t_1 - t_0 t_1^2  - 2 t_0 t_1 t_2 + t_0 + t_1 + t_2\\
				\hline 
			\end{array}$$
			
			The first five generalized Chebyshev polynomials of rank four are
			$$\begin{array}{|c|c|}
				\hline 
				P_{1, 4}(t_0,t_1,t_2,t_3) & t_0 \\
				\hline 
				P_{2, 4}(t_0,t_1,t_2,t_3) & t_0t_1-1 \\
				\hline 
				P_{3, 4}(t_0,t_1,t_2,t_3) & t_0 t_1 t_2 - t_0 - t_2 \\
				\hline 
				P_{4, 4}(t_0,t_1,t_2,t_3) & t_0 t_1 t_2 t_3 - t_0 t_1 - t_2 t_3 - t_0 t_3 + 1 \\
				\hline 
				P_{5, 4}(t_0,t_1,t_2,t_3) & t_0^2 t_1 t_2 t_3 - t_0 t_1 t_2 - t_0^2 t_1 - t_0t_2 t_3  + t_2 - t_0^2 t_3 + 2t_0\\
				\hline 
			\end{array}$$
			
			The first five generalized Chebyshev polynomials of rank $n$ for $n \geq 5$ and the first five generalized Chebyshev polynomials of infinite rank.
		\end{exmp}
	\end{subsection}
\end{section}

\begin{section}{Generalized Chebyshev polynomials and cluster algebras of type $\mathbb A$}\label{section:typeA}
	We now prove that generalized Chebyshev polynomials arise in the context of cluster algebras associated to equioriented quivers of Dynkin type $\mathbb A$. For any $r \geq 0$, we denote by $Q_r$ the quiver of Dynkin type $\mathbb A_r$ equipped with the following orientation
	$$\xymatrix{ Q_r: & 0 & \ar[l] 1 & \ar[l] 2 & \ar[l] \cdots & \ar[l] r-1.}$$
	For any $i \in [0, r-1]$, we denote by $S_i$ the simple $kQ_r$-module associated to the vertex $i$ and for any $n \in [1, r-i-1]$ we denote by $S_i^{(n)}$ the unique indecomposable module with socle $S_i$ and length $n$. By convention, $S_i^{(0)}$ denotes the zero module.
	
	We denote by $X_?$ the Caldero-Chapoton map associated to the quiver $Q_r$.
	
	\begin{prop}\label{prop:AChebyshev}
		Let $r \geq 1$, then for any $i \in [0, r-1]$ and any $n \in [1, r-i]$, we have
		$$X_{S_i^{(n)}}=P_n(X_{S_i}, \ldots, X_{S_{i+n-1}}).$$
	\end{prop}
	\begin{proof}
		We prove it by induction on $n$. If $n=1$, the result clearly holds. Let $n >1$ and fix some indecomposable $kQ_r$-module $S_i^{(n)}$. Then $i$ is not equal to $r-1$ since there are no modules of length $\geq 2$ and socle $S_{r-1}$. Thus, there is an almost split sequence of $kQ_r$-modules
		$$0 \fl S_i^{(n-1)} \fl S_i^{(n)} \oplus S_{i+1}^{(n-2)} \fl S_{i+1}^{(n)} \fl 0.$$
		It follows from \cite[Proposition 3.10]{CC} that 
		$$X_{S_i^{(n-1)}}X_{S_{i+1}^{(n-1)}}=X_{S_i^{(n)}}X_{S_{i+1}^{(n-2)}}+1.$$
		By induction, we get
		\begin{align*}
			X_{S_i^{(n)}}
				&=\frac{X_{S_i^{(n-1)}}X_{S_{i+1}^{(n-1)}}-1}{X_{S_{i+1}^{(n-2)}}}\\
				&=\frac{P_{n-1}(X_{S_i}, \ldots, X_{S_{i+n-2}})P_{n-1}(X_{S_{i+1}}, \ldots, X_{S_{i+n-1}})-1}{P_{n-2}(X_{S_{i+1}}, \ldots, X_{S_{i+n-2}})}\\
				&=P_{n}(X_{S_i}, \ldots, X_{S_{i+n-1}})
		\end{align*}
		which proves the proposition.
	\end{proof}

	\begin{corol}\label{corol:generateursA}
		Let $r \geq 1$. For every $i \in [0,r-1]$, set $u_i'=\frac{u_{i-1}+u_{i+1}}{u_i}$ with $u_{-1}=u_r=1$. Then, there are isomorphisms of $\Z$-algebras~:
		\begin{enumerate}
			\item $\mathcal A(Q_r) \simeq \Z[u_i,u_i'|i \in [0,r-1]]$;
			\item $\mathcal A(Q_r) \simeq \Z[u_0,u_i'|i \in [0,r-1]]$.
		\end{enumerate}
	\end{corol}
	\begin{proof}
		According to \cite{CC}, cluster variables in $\mathcal A(Q_r)$ are either initial cluster variables or $X_M$ where $M$ is an indecomposable $kQ_r$-module. Since every indecomposable $kQ_r$-module is of the form $S_i^{(n)}$, it follows from Proposition \ref{prop:AChebyshev} that $\mathcal A(Q_r)$ is generated as a $\Z$-algebra by the initial cluster variables $u_0, \ldots, u_{r-1}$ and the $X_{S_i}$ for $i \in [0, r-1]$. Now a direct computation shows that $X_{S_i}=u_i'$ for every $i \in [0,r-1]$. This proves the first assertion.
		
		In order to prove the second assertion, we notice that for every $i \in [0, r-2]$, $u_{i+1}=u_iu_i'-u_{i-1}$ so that $u_{i+1} \in \Z[u_j,u_j'|0 \leq j <i]$. By induction, any of the $u_i$ for $i \geq 1$ belongs to $\Z[u_0,u_i'|i \in [0,r-1]]$. This proves the second assertion.
	\end{proof}

	\begin{rmq}
		The first point of Corollary \ref{corol:generateursA} is a particular case of a result obtained independently (in a wider context) in \cite{cluster3}. Indeed, for every $i \in [0,r-1]$, the Laurent polynomial $u_i'$ is a cluster variable such that $\textbf u \setminus \ens{u_i} \sqcup \ens{u_i'}$ is the cluster obtained from the initial seed by mutating along the direction $i$. It thus follows from \cite[Corollary 1.21]{cluster3} that $u_0, \ldots, u_r, u_0', \ldots, u_r'$ generate the cluster algebra $\mathcal A(Q_r)$.
	\end{rmq}

	We now prove the algebraic independence of characters associated to simple $kQ_r$-modules. Results in this direction are usually obtained using gradings on the cluster algebra (see e.g. \cite{CK1}) but in this case, we can provide a very elementary proof using generalized Chebyshev polynomials.
	\begin{lem}
		Let $r \geq 1$ and $S_0, \ldots, S_r$ be the simple $kQ_r$-modules. Then, the set $\ens{X_{S_0}, \ldots, X_{S_{r-1}}}$ is algebraically independent over $\Z$.
	\end{lem}
	\begin{proof}
		We compute
		\begin{align*}
			P_{r+1}(X_{S_0}, \ldots, X_{S_{r-1}}, u_{r-1})
				&=u_{r-1}P_{r}(X_{S_0}, \ldots, X_{S_{r-1}})-P_{r-1}(X_{S_0}, \ldots, X_{S_{r-2}})\\
				&=X_{P_{r-1}[1]}X_{P_{r-1}}-X_{P_{r-2}}\\
		\end{align*}
		There are $k$-linear isomorphisms $\Ext^1_{\CC_{Q_r}}(P_{r-1}[1], P_{r-1}) \simeq \Ext^1_{\CC_{Q_r}}(P_{r-1}, P_{r-1}[1]) \simeq k$ and the corresponding triangles are 
		$$P_{r-1} \fl 0 \fl P_{r-1}[1] \fl P_{r-1}[1],$$
		$$P_{r-1}[1] \fl P_{r-2} \fl P_{r-1} \fl P_{r-1}[2] \simeq S_{r-1}.$$
		Applying Caldero-Keller's one-dimensional multiplication theorem \cite{CK2} gives 
		$$X_{P_{r-1}[1]}X_{P_{r-1}}=X_{P_{r-2}}+1.$$
		And thus, 
		\begin{equation}\label{eq:Pr+1}
			P_{r+1}(X_{S_0}, \ldots, X_{S_{r-1}}, u_{r-1})=1.
		\end{equation}

		We get
		$$u_{r-1}=\frac{P_{r-1}(X_{S_0}, \ldots, X_{S_{r-2}})+1}{P_r(X_{S_0}, \ldots, X_{S_{r-1}})} \in \Q(X_{S_0}, \ldots, X_{S_{r-1}}).$$
		From the second point of Corollary \ref{corol:generateursA}, we thus know that each $u_i$ belongs to $\Q(X_{S_0}, \ldots, X_{S_{r-1}})$ when $i$ runs over $[0,r-1]$. It follows that $\ens{X_{S_0}, \ldots, X_{S_{r-1}}}$ is a transcendence basis for $\Q(u_0, \ldots, u_{r-1})$.
	\end{proof}
	
	We can now give a a presentation of the cluster algebra $\mathcal A(Q_r)$ by generators and relations.
	\begin{corol}\label{corol:presentationA}
		Let $r \geq 1$ be an integer and $\mathcal A$ be a cluster algebra of Dynkin type $\A_r$. Then, there is an isomorphism of $\Z$-algebras 
		$$\mathcal A \simeq \Z[t_0, \ldots, t_{r}]/(P_{r+1}(t_0, \ldots, t_r)-1)$$
		where $t_0, \ldots, t_{r}$ are indeterminates over $\Z$.
	\end{corol}
	\begin{proof}
		Since $\mathcal A$ is of Dynkin type $\A_r$, there is an isomorphism of $\Z$-algebras $\mathcal A \simeq \mathcal A(Q_r)$. Let $\mathcal A'=\Z[t_0, \ldots, t_r]/(P_{r+1}(t_0, \ldots, t_r)-1)$. According to equation (\ref{eq:Pr+1}), there is a surjection $\mathcal A' \fl \mathcal A(Q_r)$ sending $t_r \mapsto u_{r-1}$ and $t_i \mapsto X_{S_i}$ for every $i \in [0,r-1]$.
		
		Now, consider the elements $y_0, \ldots, y_r$ in $\mathcal A'$ given by $y_0=t_0$ and 
		$$y_{i+1}=y_it_i-y_{i-1}$$
		for $i \geq 0$ with $y_{-1}=y_r=1$. In particular, for every $i \in [0,r-1]$, we have
		$$x_i=\frac{y_{i-1}+y_{i+1}}{y_i}.$$
		
		Thus,
		$$\mathcal A(Q_r)=\Z[u_i,X_{S_i}|i \in [0,r-1]]/(u_iX_{S_i}=u_{i-1}+u_{u_{i+1}})$$
		with $u_{-1}=u_r=1$ and there is an inverse map $\mathcal A(Q_r) \fl \mathcal A'$ sending $u_i \mapsto y_i$ and $X_{S_i} \mapsto x_i$.
	\end{proof}
\end{section}

\begin{section}{Generalized Chebyshev polynomials and regular modules}\label{section:regular}
	In this section, we prove that generalized Chebyshev polynomials arise in the context of characters associated to regular modules over the path algebra of a representation-infinite quiver.
	
	\begin{theorem}\label{theorem:polyreg}
		Let $Q$ be a quiver of infinite representation type. Let $\mathcal R$ be a regular component of the form $\Z\Ai/(p)$ for some $p \geq 0$ in $\Gamma(kQ)$. Let $\ens{R_i|i \in \Z/p\Z}$ denote the set of quasi-simple modules in $\mathcal R$ ordered such that $\tau R_i \simeq R_{i-1}$ for any $i \in \Z/p\Z$. Then, for any $n \geq 1$ and any $i \in \Z/p\Z$, we have
		$$X_{R_i^{(n)}}=P_n(X_{R_i}, \ldots, X_{R_{i+n-1}}).$$
		Moreover, if $p >0$, then 
		$$X_{R_i^{(n)}}=P_{n,p}(X_{R_i}, \ldots, X_{R_{i+p-1}}).$$
	\end{theorem}
	\begin{proof}
		For every $i \in \Z/p\Z$ and any $n \geq 1$, there is an almost split sequence
		$$0 \fl R_i^{(n)} \fl R_i^{(n+1)} \oplus R_{i+1}^{(n-1)} \fl R_{i+1}^{(n)} \fl 0.$$
		Thus, the corresponding characters satisfy the cluster mesh relation
		$$X_{R_i^{(n)}}X_{R_{i+1}^{(n)}}=X_{R_i^{(n+1)}}X_{R_{i+1}^{(n-1)}}+1.$$
		Now consider the epimorphism of $\Z$-algebras $\rho:\Z[x_{i}|i \in \Z] \fl \Z[X_{R_i}| i \in \Z/p\Z]$ sending $x_{i}$ to $X_{R_i}$ where the index of $R_i$ is taken in $\Z/p\Z$. It thus follows directly from the cluster mesh relation that 
		$$X_{R_i^{(n)}}=\rho(x_{i,n})=\rho(P_n(x_{i}, \ldots, x_{i+n-1})) = P_{n}(X_{R_i}, \ldots, X_{R_{i+n-1}})$$
		which proves the first assertion.
		If $p >0$, it thus follows that 
		$$P_{n}(X_{R_i}, \ldots, X_{R_{i+n-1}})=P_{n,p}(X_{R_i}, \ldots, X_{R_{i+p-1}})$$
		and the second assertion is proved.
	\end{proof}
	
	As a corollary, we get~:
	\begin{corol}\label{corol:regularinAQ}
		Let $Q$ be an acyclic quiver and $M$ be an indecomposable regular $kQ$-module with rigid quasi-composition factors. Then, $X_M \in \mathcal A(Q)$.
	\end{corol}
	\begin{proof}
		Let $n$ be the quasi-length of $M$ and $R$ be the quasi-socle of $M$. Then, the quasi-composition factors of $M$ are $R, \tau^{-1}R, \ldots, \tau^{-(n-1)}R$. Theorem \ref{theorem:polyreg} implies that
		$$X_M=P_n(X_R, X_{\tau^{-1}R}, \ldots, X_{\tau^{-(n-1)}R}).$$
		Since $\tau^{-i}R$ is rigid for every $i = 1, \ldots, n-1$, Theorem \ref{theorem:correspondenceCK2} implies that each $X_{\tau^{-i}R}$ is a cluster variable in $\mathcal A(Q)$. Thus, $X_M$ is a polynomial in cluster variables and thus belongs to $\mathcal A(Q)$.
	\end{proof}
	
	\begin{exmp}
		Fix $Q$ an affine quiver and $M$ an indecomposable regular $kQ$-module in an exceptional tube $\mathcal T$ of $\Gamma(kQ)$. It is well known that quasi-simple modules in $\mathcal T$ are always rigid (see e.g. \cite{ringel:1099}). In particular the quasi-composition factors $M$ of $N$ are rigid and thus, by Corollary \ref{corol:regularinAQ}, $X_M \in \mathcal A(Q)$.
	\end{exmp}
	
	\begin{exmp}
		Let $Q$ is a wild quiver and $M$ is a regular $kQ$-module such that $\End_{kQ}(M) \simeq k$. Assume moreover that $M$ has quasi-length at least 2. Then, it is known the quasi-composition factors of $M$ are rigid (see e.g. \cite[Chapter XVIII]{SS:volume3}) and thus, by Corollary \ref{corol:regularinAQ}, $X_M \in \mathcal A(Q)$.
	\end{exmp}

\end{section}

\begin{section}{The cluster-mesh-algebra of a stable translation quiver}\label{section:clustermeshalgebra}
	We now introduce the notion of cluster mesh algebra of a stable translation quiver. For general results concerning stable translation quiver, we refer to \cite{Riedtmann:translation}.
	
	\begin{defi}
		The \emph{cluster mesh algebra associated to the stable translation quiver $(\Gamma,\tau)$} is 
		$$\mathcal A(\Gamma,\tau)=\Z[y_i | i \in \Gamma_0]/(y_iy_{\tau(i)}=1+\prod_{j \fl i} y_j|i \in \Gamma_0)$$
		where $\ens{y_i | i \in \Gamma_0}$ is a family of indeterminates over $\Z$.
	\end{defi}
	
	\begin{exmp}
		The quiver $\Z\Ai$ is a stable translation quiver for the translation $\tau((i,n))=(i-1,n)$ for every $i \in \Z$, $n \geq 1$. Since $x_{i,n}=P_n(x_{i}, \ldots, x_{i+n-1})$ for every $i \in \Z$, $n \geq 1$, it is straightforward to see that the cluster mesh algebra of the quiver $(\Z\Ai,\tau)$ is
		$$\mathcal A(\Z\Ai,\tau)=\Z[x_{i}|i \in \Z].$$
	\end{exmp}
	
	\begin{exmp}
		Let $r \geq 1$, we consider the M\"obius band 
		$$\mathcal M_{r+1}= \Z\Ai/(\ens{(i,n)-(i+n-1,r-n+1)|i \in \Z,n \geq 1})$$
		with the translation $\tau$ induced by the translation on $\Z\Ai$. Then, we have
		$$\mathcal A(\mathcal M_{r+1},\tau)=\Z[t_0, \ldots, t_r]/(P_{r+1}(t_0,\ldots, t_r)-1).$$
		In particular, it follows from Corollary \ref{corol:presentationA} that there is a an isomorphism of $\Z$-algebras
		$$\mathcal A(Q_r) \simeq \mathcal A(\mathcal M_{r+1},\tau).$$
	\end{exmp}

	\begin{exmp}
		Let $\mathcal T_p=\Z\Ai/(p)$ be a tube of rank $p$ for some $p \geq 1$. This is a stable translation quiver for the translation induced by the translation of $\Z\Ai$. Since $x_{i,n}=P_n(x_{i}, \ldots, x_{i+n-1})$ for every $i \in \Z/p\Z$, $n \geq 1$, it is straightforward to see that the cluster mesh algebra of the tube of rank $p$ $\mathcal T_p$ is
		$$\mathcal A(\mathcal T_p,\tau)=\Z[x_{i}|i \in \Z/p\Z].$$
	\end{exmp}
	
	We now prove that if $\mathcal T$ is a tube of rank $p\geq 1$ in the Auslander-Reiten quiver $\Gamma(kQ)$ of an affine quiver and if $R_0, \ldots, R_{p-1}$ are the quasi-simple modules of $\mathcal T$, then the cluster-mesh algebra $\mathcal A(\mathcal T_p)$ is isomorphic to the $\Z$-algebra generated by $X_{R_0}, \ldots, X_{R_{p-1}}$. We will first need some technical results. 
	
	The following result is classical, we give the proof for completeness.
	\begin{lem}\label{lem:XRlinind}
		Let $Q$ be an affine quiver, $\mathcal T$ be a tube of rank $p \geq 1$ in $\Gamma(kQ)$ and $R_0, \ldots, R_{p-1}$ the quasi-simple modules in $\mathcal T$. Then $\ens{\ddim R_0, \ldots, \ddim R_{p-1}}$ is linearly independent over $\Z$.
	\end{lem}
	\begin{proof}
		We assume that the quasi-simple are ordered in such a way that $\tau R_i=R_{i-1}$ for all $i \in \Z/p\Z$. Then, for any $i \in \Z/p\Z$, we have
		$$\dim \Ext^1_{kQ}(R_i,R_j)
		=\left\{\begin{array}{rl}
			1 & \textrm{ if } j=i-1~;\\
			0 & \textrm{ otherwise.}\\
		\end{array}\right.$$
		and 
		$$\dim \Hom_{kQ}(R_i,R_j)
		=\left\{\begin{array}{rl}
			1 & \textrm{ if } j=i~;\\
			0 & \textrm{ otherwise.}\\
		\end{array}\right.$$
		
		Fix a zero linear combination $\sum_{j=0}^{p-1} \lambda_j \ddim R_j=0$ with $\lambda_0, \ldots, \lambda_{p-1} \in \Z$. For every $i \in [0,p-1]$, the linear form $\<-,\ddim R_i\>$ applied to this equality gives
		$$\lambda_i-\lambda_{i-1}=0.$$
		Thus $\lambda_i=\lambda_j$ for all $i \neq j$ and we denote by $\lambda$ this common value. Then
		$$0=\sum_{j=0}^{p-1} \lambda_j \ddim R_j=\lambda \sum_{j=0}^{p-1} \ddim R_j$$
		and $\lambda=0$.
	\end{proof}

	\begin{lem}\label{lem:XRalgind}
		Let $Q$ be an affine quiver, $\mathcal T$ be a tube of rank $p \geq 1$ in $\Gamma(kQ)$ and $R_0, \ldots, R_{p-1}$ the quasi-simple modules in $\mathcal T$. Then $\ens{X_{R_i}| i \in \Z/p\Z}$ is algebraically independent over $\Z$.
	\end{lem}
	\begin{proof}
		Let $P(t_1, \ldots, t_p) =\sum_{\nu \in \N^p} a_{\nu} t_1^{\nu_1}\cdots t_p^{\nu_p} \in \Z[t_1, \ldots, t_p]$ be a polynomial such that 
		\begin{align*}
			0	&= P(X_{R_0}, \ldots, X_{R_{p-1}})\\
				&= \sum_{\nu \in \N^p} a_{\nu} X_{\bigoplus_{i=0}^p R_i^{\nu_i}}.
		\end{align*}
		But for every $\nu \in \N^p$, $\ddim \bigoplus_{i=0}^p R_i^{\nu_i}=\sum_{i=0}^p \nu_i\ddim R_i$. Thus, Lemma \ref{lem:XRlinind} implies that $\ddim \bigoplus_{i=0}^p R_i^{\nu_i} \neq \ddim \bigoplus_{i=0}^p R_i^{\mu_i}$ is $\nu \neq \mu \in \N^p$. In particular, it follows from \cite[Proposition 4.18]{Dupont:BaseAaffine} that the $X_{\bigoplus_{i=0}^p R_i^{\nu_i}}$ are linearly independent over $\Z$ and thus, $a_{\nu}=0$ for every $\nu \in \N^p$ and $P=0$. This proves the lemma.
	\end{proof}
	
	Thus, we proved~:
	\begin{prop}
		Let $Q$ be an affine quiver and $\mathcal T$ be a tube of rank $p$ with quasi simple modules $R_0, \ldots, R_{p-1}$. Then, there is an isomorphism of $\Z$-algebras
		$$\Z[X_{R_i}|i=0,\ldots, p-1] \simeq \mathcal A(\mathcal T,\tau).$$
	\end{prop}
	\begin{proof}
		Let $\mathcal A'=\Z[X_{R_i}|i=0,\ldots, p-1]$. According to lemma \ref{lem:XRalgind}, there is an isomorphism of $\Z$-algebras $\phi:\mathcal A(\mathcal T,\tau) \fl \mathcal A'$ sending $x_{i}$ to $X_{R_i}$. Moreover, it follows from Theorem \ref{theorem:polyreg} that $\phi(x_{i,n})=X_{R_i^{(n)}}$ for every $i \in \Z/p\Z$, $n \geq 1$.
	\end{proof}

	\begin{rmq}
		If $Q$ is a wild quiver and $\mathcal R \simeq \Z\Ai$ is a regular component in $\Gamma(kQ)$ whose quasi-simple modules are denoted by $R_i, i \in \Z$, it follows from Theorem \ref{theorem:polyreg} that there is an epimorphism of $\Z$-algebras $\phi:\mathcal A(\Z\Ai,\tau) \fl \Z[R_i|i \in \Z]$ sending $x_{i}$ to $X_{R_i}$. For any $i \neq j \in \Z$, it is known that $\ddim R_i \neq \ddim R_j$ (see e.g. \cite{Zhang:compositionwild}). In particular, Caldero-Keller's denominator theorem implies that the denominator vectors of $X_{R_i}$ and $X_{R_j}$ are distinct (see \cite{CK2}). Thus, the family $\ens{X_{R_i}| i \in \Z}$ is infinite in a field $\mathcal F=\Q(u_i|i \in Q_0)$ having a finite transcendence basis over $\Q$. Thus, the $X_{R_i},i \in \Z$ are not algebraically independent over $\Z$ and $\phi$ is not an isomorphism of $\Z$-algebras.
	\end{rmq}
\end{section}

\begin{section}{Multiplications in regular components}\label{section:multiplication}
	We now prove multiplication theorems for cluster characters associated to indecomposable modules in regular components. 
	
	\begin{subsection}{The affine case}
		In this subsection, we assume that $Q$ is a quiver of affine type. Let $\mathcal T$ be a tube of rank $p$ in $\Gamma(kQ)$. As usual, we denote by $R_i, i \in \Z$ the quasi-simple modules in $\mathcal T$ ordered such that $\tau R_i \simeq R_{i-1}$ and $R_{i+p} \simeq R_i$ for every $i \in \Z$. Let $M,N$ be indecomposable regular modules contained in $\mathcal T$. Up to reordering, we can assume that $M =R_0^{(n)}$ and $N = R_j^{(m)}$ for some integers $m,n >0$ and some $j \in [0,p-1]$. 

		\begin{theorem}\label{theorem:multintubes}
			Let $Q$ be an affine quiver and $\mathcal T$ be a tube of rank $p$ in $\Gamma(kQ)$. Let $m,n>0$ be integers and $j \in [0,p-1]$. Then, for every $k \in \Z$ such that $0 < j+kp  < n$ and $m > n-j-kp$, we have the following identity~:
			$$X_{R_{j}^{(m)}}X_{R_0^{(n)}}=X_{R_0^{(m+j+kp)}}X_{R_{j}^{(n-j-kp)}}+X_{R_0^{(j+kp-1)}}X_{R_{n+1}^{(m+j+kp-n-1)}}.$$
		\end{theorem}
		\begin{proof}
			Let $k$ be in $\Z$ such that $0 < j+kp  \leq n$ and $m \geq n-j-kp$ and set $i=j+kp$. Let $Q_{i+m}$ be the quiver of Dynkin type $\A$ equipped with the following orientation
			$$\xymatrix{
				0 & \ar[l] 1 & \ar[l] \cdots & \ar[l] i+m-1
			}$$
			
			For any $r \in [0,i+m-1]$, we denote by $S_r$ (resp. $P_r$, $I_r$) the associated simple (resp. projective, injective) $kQ_{i+m}$-module. Let $X'_?$ be the Caldero-Chapoton map on $\CC_{Q_{i+m}}$. Consider the epimorphism of $\Z$-algebras:
			$$\phi:\left\{\begin{array}{rcl}
				\Z[X'_{S_r}|i \in [0,r+m-1]] & \fl & \Z[X_{R_r}|r \in [0,p-1]]\\
				X'_{S_r} & \mapsto & X_{R_{r}}
			\end{array}\right.$$
			
			Now if $S_r^{(n)}$ is any indecomposable regular $kQ_{i+m}$-module, it follows from Proposition \ref{prop:AChebyshev} that $X'_{S_r^{(n)}}=P_n(X'_{S_r}, \ldots, X'_{S_{r+n-1}})$. As Theorem \ref{theorem:polyreg} implies that $X_{R_r^{(n)}}=P_n(X_{R_r}, \ldots, X_{R_{r+n-1}})$, it follows that 
			then $\phi(X'_{S_r^{(n)}})=X_{R_{r}^{(n)}}$ for any indecomposable $kQ_{m+i}$-module $S_r^{(n)}$.
			
			Since $0 < i  < n$ and $m > n-i$, we have 
			$$\Ext^1_{kQ_{i+m}}(S_i^{(m)},S_0^{(n)}) \simeq \Hom_{kQ_{i+m}}(S_0^{(n)},\tau S_i^{(m)}) \simeq k$$
			and there is a non-split exact sequence of $kQ_{i+m}$-modules given by
			$$0 \fl S_0^{(n)} \fl S_0^{(i+m)} \oplus S_{i}^{(n-i)} \fl S_{i}^{(m)} \fl 0.$$
			Thus, applying Caldero-Keller's multiplication formula, we get
			$$X'_{S_i^{(m)}}X'_{S_0^{(n)}}=X'_{S_0^{(i+m)}}X'_{S_{i}^{(n-i)}} + X'_B$$
			where $B= \ker f \oplus \tau^{-1} \coker f$ for any $0 \neq f \in \Hom_{kQ_{i+m}}(S_0^{(n)},\tau S_i^{(m)})$. 

			Now, for every $r \in [0,m+i-1]$, we have $P_{r} \simeq S_0^{(r+1)} $ and $I_{r} \simeq S_{r}^{(m+i-r)}$ so that $\Hom_{kQ_{m+n}}(S_0^{(n)},\tau S_i^{(m)}) \simeq \Hom_{kQ_{m+n}}(P_{n-1},\tau I_i)$. Fix now some $f \neq 0$ in $\Hom_{kQ_{m+n}}(P_{n-1},\tau I_i)$. A direct computation proves that
			$\ker f \simeq P_{i-2} \simeq S_{0}^{(i-1)}$ and $\coker f \simeq \tau I_{n+1}$. Thus, 
			$$B \simeq S_0^{(i-1)} \oplus S_{n+1}^{(m+i-n-1)}$$
			and
			$$X'_{S_i^{(m)}}X'_{S_0^{(n)}}=X'_{S_0^{(i+m)}}X'_{S_{i}^{(n-i)}} + X'_{S_0^{(i-1)}}X'_{S_{n+1}^{(m+i-n-1)}}.$$
			Applying $\phi$ we get 
			$$X_{R_{i}^{(m)}}X_{R_0^{(n)}}=X_{R_0^{(m+i)}}X_{R_{j}^{(n-i)}}+X_{R_0^{(i-1)}}X_{R_{n+1}^{(m+i-n-1)}}$$
			which proves the theorem.
		\end{proof}
		
		We now prove that the multiplication provided by Theorem \ref{theorem:multintubes} can be interpreted as a Hall product in the cluster category $\CC_Q$.
		\begin{corol}\label{corol:Hallproduct}
			Let $Q$ be an affine quiver, $\mathcal T$ be a tube in $\Gamma(kQ)$ and $M,N$ be indecomposable regular modules in $\mathcal T$ such that $\Ext^1_{kQ}(N,M) \neq 0$. Then, 
			there exists two regular $kQ$-modules $B,E$ in $\mathcal T$ such that 
			$$X_MX_N=X_E+X_B$$
			and such that there exists non-split triangles 
			$$M \fl E \fl N \fl M[1] \textrm{ and } N \fl B \fl M \fl N[1]$$
			in $\CC_Q$.
		\end{corol}
		\begin{proof}
			Let $R_i, i \in \Z$ denote the quasi-simple modules in $\mathcal T$ ordered such that $\tau R_i \simeq R_{i-1}$ and $R_{i+p} \simeq R_i$ for every $i \in \Z$. Up to reordering, we can assume that $M =R_0^{(n)}$ and $N = R_j^{(m)}$ for some integers $m,n >0$ and some $j \in [0,p-1]$. 

			It is known that $\Ext^1_{kQ}(N,M) \neq 0$ if and only if there exists $k \in \Z$ such that
			$$
			\left\{\begin{array}{l}
				0 < j+kp  < n~; \\
				m > n-j-kp.
			\end{array}\right.$$
	      
			Fix any $k \in \Z$ such that $0 < j+kp < n$ and $m>n-j-kp$ and set $i=j+kp$. Since $\mathcal T$ is of rank $p$, it follows that $R_i \simeq R_j$. 
			
			We know that $\Ext^1_{kQ}(R_i^{(m)},R_0^{(n)}) \simeq \Ext^1_{kQ}(R_0^{(n)},\tau R_i^{(m)}) \neq 0$. More precisely, there is a non-split short exact sequence
			$$0 \fl M \simeq R_0^{(n)} \fl R_0^{(i+m)} \oplus R_{i}^{(n-i)} \fl R_{i}^{(m)} \simeq N \fl 0.$$
			in $kQ$-mod inducing a non-split triangle 
			$$M \fl R_0^{(i+m)} \oplus R_{j}^{(n-i)} \fl N \fl M[1]$$
			in $\CC_Q$. 
			
			Also, there is a non-zero morphism $f: R_0^{(n)} \fl \tau R_j^{(m)}$ in $kQ$-mod such that  $\ker f \simeq R_{0}^{(i-1)}$ and $\coker f \simeq R_{n}^{(m+i-n-1)}$ so that there is a non-split triangle 
			$$\tau N[-1] \fl \ker f \oplus \coker f[-1] \fl M \xrightarrow{f} \tau N$$
			in $D^b(kQ)$ inducing a non-split triangle
			$$N \fl R_{0}^{(i-1)} \oplus R_{n+1}^{(m+i-n-1)} \fl M \xrightarrow{f} N[1]$$
			in $\CC_Q$.
			
			The corollary thus follows from Theorem \ref{theorem:multintubes}.
		\end{proof}
		
		\begin{exmp}
			Let $Q$ be an affine quiver such that the Auslander-Reiten quiver $\Gamma(kQ)$ of $kQ$-mod contains a tube $\mathcal T$ of rank 4. We denote by $R_i, i \in \Z$ the quasi-simple modules in $\mathcal T$ such that $R_{i+4} \simeq R_i$ and $\tau R_i \simeq R_{i-1}$ for any $i \in \Z$.  We want to compute the product $X_{R_0}X_{R_1^{(3)}}$ where $X_?$ denotes the Caldero-Chapoton map on $\CC_Q$.
			
			Consider the quiver 
			$$\xymatrix{
				Q_8: & 0 & \ar[l] 1 & \ar[l] \cdots & \ar[l] 7
			}$$
			and denote by $X'_?$ the Caldero-Chapoton map on $\CC_{Q_8}$. Caldero-Keller multiplication theorem for $X'_?$ implies 
			$$X'_{S_0}X'_{S_1^{(3)}}=X'_{S_0^{(4)}}+X'_{S_2^{(2)}}.$$
			and 
			$$X'_{S_4}X'_{S_1^{(3)}}=X'_{S_1^{(4)}}+X'_{S_1^{(2)}}.$$
			
			Consider the epimorphism $\phi:\Z[X'_{S_i}|i=0,\ldots, 7] \fl \Z[X_{R_i}|i=0,\ldots, 3]$ sending $X'_{S_i}$ to $X_{R_i}$. Applying $\phi$ to these equalities, we obtain
			$$X_{R_0}X_{R_1^{(3)}}=X_{R_0^{(4)}}+X_{R_2^{(2)}}$$
			and 
			$$X_{R_0}X_{R_1^{(3)}}=X_{R_1^{(4)}}+X_{S_1^{(2)}}$$
			so that 
			$$2X_{R_0}X_{R_1^{(3)}}=X_{R_0^{(4)}}+X_{R_2^{(2)}}+X_{R_1^{(4)}}+X_{R_1^{(2)}}.$$
		\end{exmp}
		
		\begin{rmq}
			We now investigate the connections between Theorem \ref{theorem:multintubes} and the multiplication formula given by Xiao and Xu. Let $\mathcal T$ be a tube of rank $p \geq 1$ in the Auslander-Reiten quiver $\Gamma(kQ)$ of an affine quiver $Q$ whose quasi-simple modules are $R_i, i \in \Z$ ordered such that $R_{i+p} \simeq R_i$ and $\tau R_i \simeq R_{i-1}$ for any $i \in \Z$. We also denote by $\mathcal T$ the full subcategory of $kQ$-mod generated by objects in $\mathcal T$. Let $X_?$ be the Caldero-Chapoton map on $\CC_Q$.
			
			Let $\Aii$ be the double infinite quiver 
			$$\xymatrix{
				\Aii: \cdots & -1 \ar[l] & 0 \ar[l] & 1 \ar[l]& \ar[l] \cdots
			}$$
			and let $\rep(\Aii)$ be the category of finite dimensional representations of $\Aii$ over $k$. We denote by $S_i$ the simple representation at vertex $i$ for any $i \in \Z$ and by $S_i^{(n)}$ the representation of length $n$ with socle $S_i$ for any $i \in \Z$, $n \geq 1$. It is well known, that $\mathcal T$ is equivalent to the orbit category $\rep(\Aii)/(\tau^p)$. We choose an equivalence sending $S_i^{(n)}$ to $R_i^{(n)}$ for any $i \in \Z$ and $n \geq 1$.
			
			Fix $M$ and $N$ two indecomposable $kQ$-modules in $\mathcal T$ such that $\Ext^1_{kQ}(N,M) \neq 0$. Without loss of generality, we assume that $M \simeq R_0^{(n)}$ and $N \simeq R_j^{(m)}$ for some $m,n \geq 1$ and $j \in [0,p-1]$. We fix
			$$k_0=\max \ens{k | 0<j+kp<n \textrm{ and } m>n-j-kp}$$
			which is a well defined integer since $\Ext^1_{kQ}(N,M) \neq 0$. We set $i=j+k_0p$. Then, the wing at vertex $S_0^{(n)}$ in $\rep(\Aii)$ is equivalent to the category $kQ_{m+i}$-mod and this equivalence sends the simple representation $S_i$ of $\rep(\Aii)$ to the simple $kQ_{m+i}$-module $S_i$ for any $i \in [0,m+i-1]$.
			
			We have the following identities~:
			$$\Ext^1_{kQ}(R_j^{(m)},R_0^{(n)}) \simeq \bigoplus_{l=0}^{k_0} \Ext^1_{\Aii}(S^{(m)}_{j+lp},S_0^{(n)}) \simeq \bigoplus_{l=0}^{k_0} \Ext^1_{kQ_{m+i}}(S^{(m)}_{j+lp},S_0^{(n)}).$$
			
			Let $X'_?$ be the Caldero-Chapoton map on $\CC_{Q_{m+i}}$. For any $l \in [0,k_0]$, we have 
			$$X'_{S_{j+lp}^{(m)}}X'_{S_0^{(n)}}=X'_{S_0^{(j+lp+m)}}X'_{S_{j+lp}^{(n-j-lp)}} + X'_{S_0^{(j+lp-1)}}X'_{S_{n+1}^{(m+j+lp-n-1)}}.$$
			giving
			$$X_{R_{j}^{(m)}}X_{R_0^{(n)}}=X_{R_0^{(j+lp+m)}}X_{R_{j}^{(n-j-lp)}} + X_{R_0^{(j+lp-1)}}X_{R_{n+1}^{(m+j+lp-n-1)}}.$$
			
			Summing up these equalities when $l$ runs over $[0,k_0]$, we get
			\begin{equation}\label{eq:likeXX}
				\dim \Ext^1_{kQ}(N,M) X_MX_N=\sum_{l=0}^{k_0} \left(X_{R_0^{(j+lp+m)}}X_{R_{j}^{(n-j-lp)}} + X_{R_0^{(j+lp-1)}}X_{R_{n+1}^{(m+j+lp-n-1)}}\right).
			\end{equation}
			
			Note that all the terms occurring in the right hand side are actually middle terms of non-split triangles involving $M$ and $N$ in $\CC_Q$. In \cite{XX}, the authors give, in much more general settings, an expansion formula for $\dim \Ext^1_{kQ}(N,M) X_MX_N$ as a $\Z$-linear combination of $X_Y$ where $Y$ runs over the set of middle terms of non-split triangles involving $M$ and $N$ in $\CC_Q$. It is not clear whether their result coincides with equality (\ref{eq:likeXX}). 
			
			Note that Theorem \ref{theorem:multintubes} is slightly more precise than the multiplication provided in \cite{XX} since it allows to express $X_MX_N$ as a $\Z$-linear combination of $X_Y$'s whereas the Xiao-Xu's multiplication theorem only allows to express teh whole term $\dim \Ext^1_{kQ}(N,M) X_MX_N$ as a $\Z$-linear combination of $X_Y$'s.
		\end{rmq}
	\end{subsection}
	
	\begin{subsection}{The wild case}
		We now prove an analogue of Theorem \ref{theorem:multintubes} for regular modules over the path algebra of a wild quiver.

		\begin{theorem}\label{theorem:multinwild}
			Let $Q$ be a wild quiver and $\mathcal R$ be a regular component in $\Gamma(kQ)$. Let $R_i, i \in \Z$ denote the quasi-simple modules in $\mathcal T$ ordered such that $\tau R_i \simeq R_{i-1}$ for any $i \in \Z$. Let $m,n>0$ be integers and $j \in \Z$ such that $0 < j < n$ and $m > n-j$. Then, we have the following identity~:
			$$X_{R_{j}^{(m)}}X_{R_0^{(n)}}=X_{R_0^{(m+j)}}X_{R_{j}^{(n-j)}}+X_{R_0^{(j-1)}}X_{R_{n+1}^{(m+j-n-1)}}.$$
		\end{theorem}
		\begin{proof}
			Let $Q_{j+m}$ be the quiver of Dynkin type $\A$ equipped with the following orientation
			$$\xymatrix{
				0 & \ar[l] 1 & \ar[l] \cdots & \ar[l] j+m-1
			}$$
			and let $X'_?$ be the Caldero-Chapoton map on $\CC_{Q_{j+m}}$. 
			
			Consider the epimorphism of $\Z$-algebras:
			$$\phi:\left\{\begin{array}{rcl}
				\Z[X'_{S_r}|r \in [0,j+m-1]] & \fl & \Z[X_{R_r}|r \in [0,j+m-1]]\\
				X'_{S_r} & \mapsto & X_{R_{r}}
			\end{array}\right.$$
			
			Now if $S_r^{(n)}$ is any indecomposable regular $kQ_{j+m}$-module, we know from cluster mesh relations that 
			$\phi(X'_{S_r^{(n)}})=X_{R_{r}^{(n)}}$.
			
			Since $0 < j  < n$ and $m >n-j$, we have isomorphisms of $k$-vector spaces $\Ext^1_{Q_{j+m}}(S_j^{(m)},S_0^{(n)}) \simeq \Hom_{Q_{i+m}}(S_0^{(n)},\tau S_j^{(m)}) \simeq k$.
			
			As in the proof of Theorem \ref{theorem:multintubes}, we get
			$$X'_{S_j^{(m)}}X'_{S_0^{(n)}}=X'_{S_0^{(j+m)}}X'_{S_{j}^{(n-j)}} + X'_{S_0^{(j-1)}}X'_{S_{n+1}^{(m+j-n-1)}}.$$
			and thus, applying $\phi$, we obtain
			$$X_{R_{j}^{(m)}}X_{R_0^{(n)}}=X_{R_0^{(m+j)}}X_{R_{j}^{(n-j)}}+X_{R_0^{(j-1)}}X_{R_{n+1}^{(m+j-n-1)}}$$
			which establishes the theorem.
		\end{proof}
	\end{subsection}
\end{section}

\section*{Acknowledgements}
	The author would like to thank P. Caldero and I. Reiten for corrections and remarks about this paper. He would also like to thank the coordinators of the Liegrits network for organizing his stay at the NTNU in Trondheim where this work was done. The author is also very grateful to anonymous referees for their very interesting comments and suggestions.


\newcommand{\etalchar}[1]{$^{#1}$}

\end{document}